\documentclass[preprint,review,10pt]{elsarticle}
\usepackage{amsfonts}
\usepackage{amsmath,amssymb}
\usepackage{graphicx}
\usepackage{setspace}
\usepackage[left=2cm,top=2cm,right=2cm]{geometry}
\usepackage{subfig}
\newtheorem{lemma}{Lemma}

\newtheorem{theorem}{Theorem}
\numberwithin{equation}{section} \journal{XYZ}
\begin{document}
\title{Some Approximation Properties of the Generalized Baskakov operators}
\author[label1,label2,slabel*]{Prashantkumar Patel}
\ead{prashant225@gmail.com}
\author[label1,label3]{Vishnu Narayan Mishra}
\ead{vishnu\_narayanmishra@yahoo.co.in;
vishnunarayanmishra@gmail.com}
\author[label4]{Mediha \"{O}rkc\"{u}}
\ead{medihaakcay@gazi.edu.tr}
\address[label1]{Department of Applied Mathematics \& Humanities,
S. V. National Institute of Technology, Surat-395 007 (Gujarat),
India}
\address[label2]{Department of Mathematics, St. Xavier's College(Autonomous), Ahmedabad-380 009 (Gujarat), India}
\address[label4]{Department of Mathematics, Faculty of Sciences, Gazi University, Teknikokullar, 06500 Ankara, Turkey}
\address[label3]{L. 1627 Awadh Puri Colony Beniganj, Phase-III, Opposite - Industrial Training Institute (ITI), Ayodhya Main Road, Faizabad, Uttar Pradesh 224 001, India}
\fntext[label*]{Corresponding authors}
\begin{abstract}
The present paper deals with a generalization of the Baskakov
operators. Some direct theorems, asymptotic formula and
$A$-statistical convergence are established. Our results are based
on a $\rho$ function. These results include the preservation
properties of the classical Baskakov operators.
\end{abstract}
 \maketitle
\section{Introduction}\label{s:1}
In \cite{baskakov1957instance}, Baskakov discussed the following
positive linear operators on the unbounded interval $[0,\infty)$,
\begin{equation}\label{m17.eq1}
V_n(f;x) = \sum_{k=0}^{\infty} v_{k,n}(x)
f\left(\frac{k}{n}\right),~~~~ x\in [0,\infty), n\in \mathbf{N},
\end{equation}
where $f$ is an appropriate function defined on the unbounded
interval $[0,\infty)$, for which the above
series is convergent and $\displaystyle v_{n,k}(x) = {n+k-1 \choose k} \frac{x^k}{(1+x)^{n+k}}$.  \\
In 2011, C{\'a}rdenas-Morales \textit{et al.}
\cite{cardenas2011bernstein} introduced a generalized Bernstein
operators by fixing $e_0$ and $e_1$, given by
\begin{equation}\label{m17.eq2}
L_n(f;x) = \sum_{k=0}^{n} { n \choose k} x^{2k}(1-x^2)^{n-k}
f\left(\sqrt{\frac{k}{n}}\right),~~~~ x\in [0,1], n\in \mathbf{N},
\end{equation}
where $f\in C[0,1]$. This is a special case of the operators
$B_{n}^{\tau} f= B_n\left(f\circ\tau^{-1}\right)\circ \tau$, for
$\tau=e_2$,
 where $B_n$ is the classical Bernstein operators.\\
Consider a real valued function $\rho$ on $[0,\infty)$ satisfied
following two conditions:
\begin{enumerate}
 \item\label{m17.item1} $\rho$ is a continuously differentiable function on $[0,\infty)$,
 \item\label{m17.item2} $\rho(0)=0$ and $\displaystyle \inf_{x\in [0,\infty)} \rho'(x)\geq 1$.
 \end{enumerate}
 Throughout the manuscript, we denote the above two conditions as  $c_1$ and $c_2$.
Recently, In \cite{aral2014generalized} the following
generalization of Sz{\'a}sz--Mirakyan operators are constructed,
\begin{equation}\label{m17.eq3}
S_n^{\rho}(f;x) = \exp\left(-n\rho(x)\right) \sum_{k=0}^{\infty}
\left(f\circ\rho^{-1}\right)\left(\frac{k}{n}\right)\frac{\left(n\rho(x)\right)^k}{k!},~~~~
x\in [0,1], n\in \mathbf{N}.
\end{equation}
Notice that if $\rho=e_1$ then the operators (\ref{m17.eq3})
reduces to the well known Sz{\'a}sz--Mirakyan operators.
 Aral \textit{et al.} \cite{aral2014generalized}, gave quantitative type theorems in order to obtain the degree of
 weighted convergence with the help of a weighted modulus of continuity constructed using the function $\rho$
 of the operators (\ref{m17.eq3}).
  Very recently, some researchers have discussed shape preserving properties of the generalized Bernstein,
  Baskakov and Sz\'{a}sz-Mirakjan operators in
  \cite{zhang2009preservation,Wang2014Shape,Wang2014shapep,mishra2013statistical}.
\section{Construction of the Operators}
This motivated us to generalize the Baskakov operators
(\ref{m17.eq1}) as
\begin{eqnarray}\label{m17.eq4}
V_n^{\rho} (f;x) &=& \sum_{k=0}^{\infty}
\left(f\circ\rho^{-1}\right)\left(\frac{k}{n}\right)
                {n+k-1 \choose k} \frac{\left(\rho(x)\right)^k}{\left(1+\rho(x)\right)^{n+k}}\\
&=& \left(V_n^{\rho} \left(\left(f\circ\rho^{-1}\right)\circ\rho\right)\right)(x)\nonumber\\
&=& \sum_{k=0}^{\infty}
f\left(\rho^{-1}\left(\frac{k}{n}\right)\right)
v_{\rho,n,k}(x),\nonumber
\end{eqnarray}
where $n\in \mathbf{N}$, $x\in [0,\infty)$ and $\rho$ is a function defined as in conditions $c_{\ref{m17.item1}}$ and $c_{\ref{m17.item2}}$.\\
Observe that, $V_n^{\rho}(f;x) = V_n(f;x)$ if $\rho=e_1$. In fact,
direct calculation gives that
\begin{equation}\label{m17.moment1} \displaystyle V_n^{\rho}(e_0;x) = 1;\end{equation}
\begin{equation}\label{m17.moment2}\displaystyle V_n^{\rho}(\rho;x) = \rho(x);\end{equation}
\begin{equation}\label{m17.moment3}\displaystyle V_n^{\rho}(\rho^2;x) = \rho^2(x)+ \frac{\rho^2(x) + \rho(x)}{n}.\end{equation}
In this manuscript, we are dealing with approximation properties
the operators (\ref{m17.eq4}). In the next section, we establish
some direct results using generalized modulus of continuity. The
Voronovskaya Asymptotic formula and $A$-Statistical convergence of
the operators $V^{\rho}_n$ are discuss in Section 4 and 5.
\section{Direct Theorems}
Consider $\phi^2(x) =1+\rho^2(x)$. Notice that, $\displaystyle
\lim_{x\to \infty} \rho(x) =\infty$ because the condition
$c_{\ref{m17.item2}}$. Denote $B_{\phi}\left([0,\infty)\right)$ as
set of all real valued function on $[0,\infty)$ such that
$\displaystyle |f(x)|\leq M_f\phi(x)$, for all $x\in [0,\infty)$,
where $M_f$ is a constant depending on $f$. Observe that,
$B_{\phi}\left([0,\infty)\right)$ is norm linear space with the
norm $\displaystyle \|f\|_{\phi}
=\sup\left\{\frac{|f(x)|}{\phi(x)}: x \in [0,\infty)\right\}$.
Also, $C_{\phi}\left([0,\infty)\right)$ as the set all continuous
function in $B_{\phi}\left([0,\infty)\right)$ and
$$C_{\phi}^k \left([0,\infty)\right)=\{f\in C_{\phi}\left([0,\infty)\right):
 \lim_{x\to \infty}\frac{f(x)}{\phi(x)}=k_f, k_f \textrm{ is constant depaneding on } f\}.$$
Let $U_{\phi} \left([0,\infty)\right)$ be the space of functions
$f\in C_{\phi} \left([0,\infty)\right)$, such that $\displaystyle
\frac{f(x)}{\phi(x)}$ is uniformly continuous. Also,
$\displaystyle C_{\phi}^k \left([0,\infty)\right)\subset U_{\phi}
\left([0,\infty)\right) \subset C_{\phi}
\left([0,\infty)\right)\subset B_{\phi} \left([0,\infty)\right)$.
\begin{lemma}\cite{Gadziev1974conve}\label{m17.lemma1}
The positive linear operators $L_n: C_{\phi}
\left([0,\infty)\right)\to B_{\phi} \left([0,\infty)\right)$ for
all $n\geq 1$ if and only if the inequality
$$|L_n (\phi; x)| \leq  K_n\phi(x),~~~~ x \in [0,\infty), ~n \geq 1,$$
holds, where $K_n$ is a positive constant.
\end{lemma}
\begin{theorem}\cite{Gadziev1974conve}\label{m17.theorem1}~
Let the sequence of linear positive operators $(L_n)_{n\geq1}$,
$L_n: C_{\phi} \left([0,\infty)\right)\to B_{\phi}
\left([0,\infty)\right)$
  satisfy the three conditions
$$\lim_{n\to \infty} \|L_n\rho^i-\rho^i\|=0,~~ i=0,1,2.$$
Then
$$\lim_{n\to \infty} \|L_nf-f\|=0,$$
for any $f\in  C_{\phi}^k \left([0,\infty)\right)$.
\end{theorem}
By (\ref{m17.moment1}), (\ref{m17.moment3}) and Lemma
\ref{m17.lemma1}, $V_n^{\rho}$ is linear positive operators from
 $C_{\phi} \left([0,\infty)\right)$ to $B_{\phi} \left([0,\infty)\right)$.
\begin{theorem}
For each function $f\in C_{\phi}^k \left([0,\infty)\right)$, we
have
\begin{equation}\label{m17.eq7}\lim_{n\to \infty} \|V_n^{\rho}(f;\cdot) - f \|_{\phi} =0.\end{equation}
\end{theorem}
\textbf{Proof: } From (\ref{m17.moment1}) and (\ref{m17.moment2}),
we write
$$ \|V_n^{\rho}(1;\cdot)-1\|_{\phi}= 0 \textrm { and } \|V_n^{\rho}(\rho;\cdot)-\rho\|_{\phi}= 0.$$
Also,
\begin{equation}\label{m17.eq8}\|V_n^{\rho}(\rho^2;\cdot)-\rho^2\|_{\phi}
             = \sup_{x\in [0,\infty)}\frac{\rho^2(x)+\rho(x)}{n(1+\rho^2(x))}\leq \frac{2}{n}~. \end{equation}
Therefore, we have $$\|V_n^{\rho}(\rho^i;\cdot)-\rho^i\|_{\phi}\to
0 \textrm{ as } n\to \infty, \textrm{ for } i=0,1,2.$$
Hence by Theorem \ref{m17.theorem1}, the equation (\ref{m17.eq7}) is also true.\\
In \cite{Holhos2008quatitati} the following weighted modulus of
continuity is defined
\begin{equation}\omega_\rho(f;\delta):= \omega(f;\delta)_{[0,\infty)}
        =\sup_{\substack{x,t\in [0,\infty) \\|\rho(x)-\rho(t)|\leq \delta}}\frac{|f(t)-f(x)|}{|\phi(t)-\phi(x)|}\end{equation}
for each $f\in C_{\phi}\left([0,\infty)\right)$ and for every $\delta>0$.\\
We call the function $\omega_{\rho}(f; \delta)$ weighted modulus
of continuity. We observe that $\omega_{\rho}(f; 0) = 0$ for every
$f\in C_{\phi}\left([0,\infty)\right) $ and that the function
$\omega_{\rho}(f;\delta)$ is nonnegative and nondecreasing with
respect to $\delta$ for $f \in C_{\phi}\left([0,\infty)\right)$.
Here, we consider the spaces $C_{\phi}^k\left([0,\infty)\right)$,
$U_{\phi}\left([0,\infty)\right)$,
$C_{\phi}\left([0,\infty)\right)$ and
$B_{\phi}\left([0,\infty)\right)$ having the conditions
$c_{\ref{m17.item1}}$ and $c_{\ref{m17.item2}}$. Under these
conditions, $|x -t| \leq |\rho(x)-\rho(t)|$, for every $x, t\in
[0,\infty)$ is true.
\begin{lemma}\cite{Holhos2008quatitati}\label{m17.lemma2}
$\displaystyle \lim_{\delta\to 0}\omega_{\rho}(f;\delta) = 0,$ for
every $f\in U_{\phi}\left([0,\infty)\right)$.
\end{lemma}
\begin{theorem}\cite{Holhos2008quatitati}\label{m17.theorem2}
Let $L_n : C_{\phi}\left([0,\infty)\right) \to
B_{\phi}\left([0,\infty)\right)$be a sequence of positive linear
operators with
\begin{eqnarray}
\|L_n(\rho^0)-\rho^0\|_{\phi^0} &=&a_n,\label{m17.eq9}\\
\|L_n(\rho)-\rho\|_{\phi^{\frac{1}{2}}} &=&b_n,\label{m17.eq10}\\
\|L_n(\rho^2)-\rho^2\|_{\phi} &=&c_n,\label{m17.eq11}\\
\|L_n(\rho^3)-\rho^3\|_{\phi^{\frac{3}{2}}}&=&d_n\label{m17.eq12},
\end{eqnarray}
where $a_n$, $b_n$, $c_n$ and $d_n$ tends to zero as $n\to
\infty$. Then
$$\|L_n(f) -f \|_{\phi^{\frac{3}{2}}}= (7+4a_n+2c_n) \omega_{\rho}(f;\delta_n)+ a_n\|f\|_{\phi}, $$
for all $f\in  C_{\phi}\left([0,\infty)\right)$, where
$$\delta_n = 2\sqrt{(a_n + 2b_n + c_n) (1+a_n)} + a_n + 3b_n + 3c_n + d_n.$$
\end{theorem}
\begin{theorem}\label{m17.theorem3}
For all $f\in  C_{\phi}\left([0,\infty)\right)$, we have
\begin{equation*}
\|V_n^{\rho}(f;\cdot)-f\|_{\phi^{\frac{3}{2}}}\leq
\left(7+\frac{4}{n}\right)\omega_{\rho}\left(f,\frac{2\sqrt{2}}{\sqrt{n}}+
\frac{16}{n}\right).
\end{equation*}
\end{theorem}
\textbf{Proof: } Notice that
\begin{equation}\|V_n^{\rho}(\rho^0;\cdot)-\rho^0\|_{\phi^0} = 0 =a_n\end{equation}
and
\begin{equation}\|V_n^{\rho}(\rho;\cdot)-\rho\|_{\phi^{\frac{1}{2}}} = 0 =b_n.\end{equation}
Form equation (\ref{m17.eq8}), we have
\begin{equation}c_n=\|V_n^{\rho}(\rho^2;\cdot)-\rho^2\|_{\phi}\leq \frac{2}{n}.\end{equation}
Now,
\begin{equation}
V_n^{\rho}(\rho^3;x)=\frac{\rho(x)}{n^2}+\frac{3 (1+n)
\rho(x)^2}{n^2}+\frac{\left(2+3 n+n^2\right) \rho(x)^3}{n^2}.
\end{equation}
We can write
\begin{eqnarray*}
d_n= \|V_n^{\rho}(\rho^3;\cdot)-\rho^3\|_{\phi^{\frac{3}{2}}}&=&
        \sup_{x\in[0,\infty)}\left\{\frac{\rho(x)}{n^2\left(1+\rho^2(x)\right)^{\frac{3}{2}}}\right.\\
        &&\left.+\frac{3 (1+n) \rho(x)^2}{n^2\left(1+\rho^2(x)\right)^{\frac{3}{2}}}+\frac{(2+3 n) \rho(x)^3}{n^2\left(1+\rho^2(x)\right)^{\frac{3}{2}}}\right\}\\
&\leq& \frac{1}{n}+\frac{4}{n}+ \frac{5}{n}= \frac{10}{n}.
\end{eqnarray*}
Observe that, the condition (\ref{m17.eq9})-(\ref{m17.eq12}) are
satisfied, therefore by theorem \ref{m17.theorem2}, we have
\begin{equation*}
\|V_n^{\rho}(f;\cdot)-f\|_{\phi^{\frac{3}{2}}}\leq
\left(7+\frac{4}{n}\right)\omega_{\rho}\left(f,\frac{2\sqrt{2}}{\sqrt{n}}+
\frac{16}{n}\right).
\end{equation*}
This completes the proof of Theorem \ref{m17.theorem3}.
\begin{theorem}
For all $f\in  U_{\phi}^k \left([0,\infty)\right)$, we have
$\displaystyle \lim_{n\to
\infty}\|V_n^{\rho}(f;\cdot)-f\|_{\phi^{\frac{3}{2}}}=0. $
\end{theorem}
The proof follows from the Theorem \ref{m17.theorem3} and Lemma
\ref{m17.lemma2}.
\section{Voronovskaya Asymptotic formula}
Now we give the following Voronovskaya type theorem. We use the
technique developed in
\cite{cardenas2011bernstein,aral2014generalized}.
\begin{theorem}
Let $f \in C[0,\infty)$, $x \in [0,\infty)$ and suppose that the
first and second derivatives of $f\circ \rho^{-1}$ exist at
$\rho(x)$. If the second derivative of $f\circ \rho^{-1}$ is
bounded on $[0,\infty)$ then we have
\begin{equation}\lim_{n\to \infty}n (V_n^{\rho}(f;x)-f(x))= \frac{1}{2}\rho(x) (1 + \rho(x)) \left(f\circ \rho^{-1}\right)''\left(\rho(x)\right).\end{equation}
\end{theorem}
\textbf{Proof: } By the Taylor expansion of $f\circ\rho^{-1}$ at
the point $\rho(x) \in [0,\infty)$, there exists $\xi$ lying
between $x$ and $t$ such that
\begin{eqnarray}\label{m17.eq14}
f(t)&=& \left(f\circ\rho^{-1}\right)\left(\rho(t)\right)=\left(f\circ\rho^{-1}\right)\left(\rho(x)\right)+ \left(f\circ\rho^{-1}\right)'\left(\rho(x)\right)\left(\rho(t)-\rho(x)\right)\nonumber\\
&&+\frac{1}{2}\left(f\circ
\rho^{-1}\right)''\left(\rho(x)\right)\left(\rho(t)-\rho(x)\right)^2+
\lambda_x(t) \left(\rho(t)-\rho(x)\right)^2,
\end{eqnarray}
where
\begin{equation}\label{m17.eq13}
\lambda_x(t) =
\frac{\left(f\circ\rho^{-1}\right)'\left(\rho(\xi)\right)-\left(f\circ\rho^{-1}\right)''\left(\rho(x)\right)}{2}.
\end{equation}
Note that, the assumptions on $f$ together with definition (\ref{m17.eq13}) ensure that $|\lambda_x(t)|\leq  M$ for all $t$ and converges to zero as $t\to x$.\\
Applying the operators (\ref{m17.eq4}) to the above equation
(\ref{m17.eq14}) equality, we get
\begin{eqnarray}\label{m17.eq15}
 V_n^{\rho}(f;x)-f(x) &=& \left(f\circ\rho^{-1}\right)'\left(\rho(x)\right)V_n^{\rho}\left(\left(\rho(t)-\rho(x)\right);x\right)\nonumber\\
&&+\frac{1}{2}\left(f\circ
\rho^{-1}\right)''\left(\rho(x)\right)V_n^{\rho}\left(\left(\rho(t)-\rho(x)\right)^2;x\right)\notag\\
&&+ V_n^{\rho}\left(\lambda_x(t)
\left(\rho(t)-\rho(x)\right)^2;x\right),
\end{eqnarray}
by equations (\ref{m17.moment1}), (\ref{m17.moment2}) and
(\ref{m17.moment3}), we have
$$\lim_{n\to \infty} n V_n^{\rho}\left(\left(\rho(t)-\rho(x)\right);x\right) = 0;$$
$$\lim_{n\to \infty} n V_n^{\rho}\left(\left(\rho(t)-\rho(x)\right)^2;x\right) = \rho(x) (1 + \rho(x)).$$
Therefore,
\begin{eqnarray}\label{m17.eq16}
 \lim_{n\to \infty}n\left(V_n^{\rho}(f;x)-f(x)\right) &=& \frac{1}{2}\rho(x) (1 + \rho(x)) \left(f\circ \rho^{-1}\right)''\left(\rho(x)\right)\nonumber\\
 &&+ \lim_{n\to \infty} n\left(V_n^{\rho}\left(\lambda_x(t)
 \left(\rho(t)-\rho(x)\right)^2;x\right)\right).
\end{eqnarray}
Now we estimate the last term on the right hand side of the above
equality. Let $\epsilon > 0$ and choose $\delta > 0$ such that
$|\lambda_x (t)| < \epsilon$ for $|t- x| <\delta$. Also it is
easily seen that by condition $c_{\ref{m17.item2}}$,
$|\rho(t)-\rho(x)| = \rho(\eta) |t - x| \geq |t- x|$. Therefore,
if $|\rho(t)-\rho(x)| < \delta$, then $|\lambda_x(t)
(\rho(t)-\rho(x))^2|< \epsilon \left(\rho(t)-\rho(x)\right)^2$ ,
while if $|\rho(t) -\rho(x)|\geq \delta$, then since
$|\lambda_x(t)|\leq M$ we have $
|\lambda_x(t)\left(\rho(t)-\rho(x)\right)^2|\leq
\frac{M}{\delta^2}\left(\rho(t)- \rho(x)\right)^4$. So we can
write
\begin{eqnarray}\label{m17.eq17}
V_n^{\rho}\left(\lambda_x(t)
\left(\rho(t)-\rho(x)\right)^2;x\right)&<&\epsilon
\left(V_n^{\rho}\left(
\left(\rho(t)-\rho(x)\right)^2;x\right)\right)\notag\\
&&+\frac{M}{\delta^2} \left(V_n^{\rho}\left(
\left(\rho(t)-\rho(x)\right)^4;x\right)\right).
\end{eqnarray}
Direct calculations show that
$$V_n^{\rho}\left( \left(\rho(t)-\rho(x)\right)^4;x\right)=O\left(\frac{1}{n^2}\right).$$
Hence $$\lim_{n\to \infty} n V_n^{\rho}\left(\lambda_x(t)
\left(\rho(t)-\rho(x)\right)^2;x\right)=0,$$ which completes the
proof of the Theorem 8.
\section{$A$-Statistical Convergence}
Now, we introduces some notation and the basic definitions, which
used in this section. Let $A=(a_{ij})$ be an infinite summability
matrix. For given sequence $x=(x_n)$, the $A-$transform to $x$,
denoted by $Ax=\left(\left(Ax\right)_j\right)$, is given by
$\displaystyle \left(Ax\right)_j=\sum_{n=1}^{\infty}a_{jn}x_n,$
provided the series
converges for each $j$. We say that $A$ is regular, if $\displaystyle\lim_j\left(Ax\right)_j=L$ whenever $\displaystyle\lim_j x_j =L$ \cite{Hardy1949diver}.\\
Now, we assume that $A$ is a nonnegative regular summability
matrix and $K$ is a subset of $\mathbf{N}$, the set of all natural
numbers. The $A$-density of $K$ is defined by $\displaystyle
\delta_A(K) = \lim_j \sum_{n=1}^{\infty} a_{jn} \chi_K(n)$
provided the limit exists, where $\chi_K$ is the characteristic
function of $K$. Then the sequence $x = (x_n)$ is said to be
$A$-statistically convergent to the number $L$ if, for every
$\epsilon> 0$, $\delta_A\{n \in \mathbf{N}  : |x_n - L| \geq
\epsilon \} = 0$; or equivalently $\displaystyle \lim_j
\sum_{n:|x_n-L|\geq \epsilon}a_{jn}=0$.
We denote this limit by $st_A -lim~ x = L$ \cite{freedman1981densities,connor1989strong,connor1996statistical,miller1995measure}.\\
For the case in which $A = C_1$, the Ces\`{a}ro matrix,
$A$-statistical convergence reduces to statistical convergence
\cite{fast1951convergence,fridy1985statistical,fridy1997statistical}.
Also, taking $A = I$, the identity matrix, $A$-statistical
convergence coincides with the ordinary convergence. We also note
that if $A = (a_{jn})$ is a nonnegative regular summability matrix
for which $\displaystyle \lim_j \max_n\{a_{jn}\} = 0$, then
$A$-statistical convergence is stronger than convergence
\cite{kolk1993matrix}. A sequence $x = (x_n)$ is said to be
$A$-statistically bounded provided that there exists a positive
number $M$ such that $\delta_A\{n\in \mathbf{N} : |x_n| \leq  M\}
= 1$. Recall that $x = (x_n)$ is A-statistically convergent to L
if and only if there exists a subsequence ${x_{n(k)}}$ of $x$ such
that $\delta_A\{n(k) : k \in\mathbf{N}\} = 1$ and $\displaystyle
\lim_k x_{n(k)} = L$ (see
\cite{miller1995measure,kolk1993matrix}).
Note that, the concept of $A$-statistical convergence is also given in normed spaces \cite{kolk1991statistical}.\\
In this section, we denote $B_{\phi}\left([0,\infty)\right)$ by
$B_{\phi}$ and $C_{\phi}\left([0,\infty)\right)$ by $C_{\phi}$.
Assume $\phi_1(x)$ and $\phi_2(x)$ be weight functions satisfying
$\displaystyle \lim_{|x|\to \infty} \frac{\phi_1(x)}{\phi_2(x)}=
0.$ If $T$ is a positive operators such that $T:C_{\phi_1} \to
B_{\phi_2}$,
then the operators norm $\|T\|_{C_{\phi_1} \to B_{\phi_2}}$ is given by $\|T\|_{C_{\phi_1} \to B_{\phi_2}}:=\sup_{\|f\|_{\phi_1=1}}\|T f \|_{\phi_2}$.\\
\begin{theorem}\cite[Thorem 6]{Duman2004stat}\label{m17.theorem4}  Let $A = (a_{jn})$ be a non-negative regular
 summability matrix, let $\{T_n\}$ be a sequence of positive linear operators from $C_{\phi_1}$ into $B_{\phi_2}$
 and  assume that $\phi_1(x)$ and $\phi_2(x)$ be weight functions satisfying $\displaystyle \lim_{|x|\to \infty} \frac{\phi_1(x)}{\phi_2(x)}= 0.$ Then
\begin{equation}
st_A-\lim_n\|T_n f-f\|_{\phi_2}=0, \textrm{ for all } f\in
C_{\phi_1}
\end{equation}
if and only if
\begin{equation}
st_A-\lim_n\|T_n \rho^{v}-\rho^{v}\|_{\phi_1}=0,~~~ v=0,1,2.
\end{equation}
\end{theorem}
With the help of Theorem \ref{m17.theorem4} we write the following
Korovkin type theorem.
\begin{theorem}\label{m17.theorem10}
Let $A = (a_{jn})$ be a non-negative regular summability matrix,
let $\{V_n\}$ be a sequence of positive linear operators from
$C_{\phi_1}$ into $B_{\phi_2}$ as defined in (\ref{m17.eq4}) and
assume that $\phi_1(x)$ and $\phi_2(x)$ be weight functions
satisfying $\displaystyle \lim_{|x|\to \infty}
\frac{\phi_1(x)}{\phi_2(x)}= 0.$ Then
\begin{equation}\label{m17.eq18}
st_A-\lim_n\|V_n(f,\cdot)-f\|_{\phi_2}=0, \textrm{ for all } f\in
C_{\phi_1}.
\end{equation}
%if and only if
\end{theorem}
\textbf{Proof: } By theorem \ref{m17.theorem4} it is sufficient to
prove that,
\begin{equation}\label{m17.eq19}
st_A-\lim_n\|V_n(\rho^{v},\cdot)-\rho^{v}\|_{\phi_1}=0,~~~
v=0,1,2.
\end{equation}
It clear that
$$\|V_n(\rho^0,\cdot)-\rho^{0}\|_{\phi_1}=0 \textrm{ and } \|V_n(\rho,\cdot)-\rho\|_{\phi_1} =0.$$
Hence, equation (\ref{m17.eq19}) is true for $v=0,1$.\\
Now, for $v=2$
\begin{equation}
\|V_n(\rho^2,\cdot)-\rho^2\|_{\phi_1} \leq \frac{2}{n}.
\end{equation}
Due to the equality $\displaystyle st_A-\lim_n\frac{1}{n}=0$, the
above inequality implies that
\begin{equation}
st_A-\lim_n\|V_n(\rho^{2},\cdot)-\rho^{2}\|_{\phi_1}=0,
\end{equation}
which completes the proof the Theorem \ref{m17.theorem10}.
\section*{Conclusions}
We constructed sequences of the Baskakov operators which are based
on a function $\rho$. This function not only characterizes the
operators but also characterizes the Korovkin set $\{1, \rho,
\rho^2\}$ in a weighted function space.\\
Our results include the following: The rate of convergence of
these operators to the identity operator on weighted spaces which
are constructed using the function $\rho$ and which are subspaces
of the space of continuous functions on $[0,\infty)$. We gave
quantitative type theorems in order to obtain the degree of
weighted convergence with the help of a weighted modulus of
continuity constructed using the function $\rho$ and the study of
$A$-statistical convergence of the sequence.
\newcommand{\nosort}[1]{}

\biboptions{numbers,sort&compress}
%\bibliographystyle{plain}  % Style BST file
%  \bibliography{Patel,Mishra and Orkcu}

\end{document}